# On the performances of a new thresholding procedure using tree structure

## Florent Autin

*Université Aix-Marseille 1, CNRS UMR 6632,*
*C.M.I., 39 rue F. Joliot Curie, 13453 Marseille Cedex 13.*
*e-mail:* autin@cmi.univ-mrs.fr

**Abstract:** This paper deals with the problem of function estimation. Using the white noise model setting, we provide a method to construct a new wavelet procedure based on thresholding rules which takes advantage of the dyadic structure of the wavelet decomposition. We prove that this new procedure performs very well since, on the one hand, it is adaptive and near-minimax over a large class of Besov spaces and, on the other hand, the maximal functional space (maxiset) where this procedure attains a given rate of convergence is very large. More than this, by studying the shape of its maxiset, we prove that the new procedure outperforms the hard thresholding procedure.

**AMS 2000 subject classifications:** 62G05, 62G07.
**Keywords and phrases:** Besov spaces, estimation, maxiset, minimax risk, rate of convergence, thresholding methods, tree structure.

Received March 2008.

## 1. Introduction

In nonparametric statistics, a lot of statisticians are interested with the estimation of a function from noisy observations. In this setting, people look for data-driven procedures able to perform very well, that is to say, for procedures very close to the target function. To reach this goal, a criterion is necessary to measure the performance of any procedure. One of the most usual way to measure this performance is to evaluate its maximum risk over a functional space $\mathcal{F}$ which the unknown signal is supposed to belong. In the $\mathbb{L}_2$-case, the maximum risk of any procedure of estimation $\hat{f}$ on $\mathcal{F}$ is the quantity

$$\mathcal{R}_{\mathcal{F}}(\hat{f}, \epsilon) := \sup_{f \in \mathcal{F}} \mathbb{E} \|\hat{f} - f\|_2^2,$$

where $\epsilon > 0$ is the noise level. In the minimax setting, the main goal is to provide procedures which are as close as possible to the $\mathcal{F}$-minimax rate $\rho_{\mathcal{F}}$ defined for any $\epsilon > 0$ by

$$\rho_{\mathcal{F}}(\epsilon) := \inf_{\hat{f}} \mathcal{R}_{\mathcal{F}}(\hat{f}, \epsilon) = \inf_{\hat{f}} \sup_{f \in \mathcal{F}} \mathbb{E} \|\hat{f} - f\|_2^2,$$

where the infimum is taken over all the data driven procedures. The minimax theory has been largely developed since the 1980-ies. A lot of minimax results have been obtained for Sobolev classes, Hölder classes and Besov classes.

Nevertheless, it appears that the minimax approach is not realistic since it requires the statistician to know the functional space $\mathcal{F}$ containing the unknown target function. Hence this point of view seems quite subjective and debatable. Moreover building an estimator adapted to the *worst* functions of $\mathcal{F}$ is not what applied statisticians are especially interested with.

Keeping in mind these minimax drawbacks, Cohen, De Vore, Kerkyacharian and Picard [6] have suggested an alternative approach to measure the performance of an estimation procedure: the maxiset point of view consists in exhibiting the largest subspace of $\mathbb{L}_2$ (maxiset) over which an estimator attains a given rate of convergence. To prove that a functional space $\mathcal{A}$ is the maxiset of a chosen procedure for a rate $r = (r_\epsilon)_\epsilon$ requires two steps. The first step is to prove that

$$\sup_{\epsilon > 0} r_\epsilon^{-1} \mathbb{E}\|\hat{f} - f\|_2^2 < \infty \implies f \in \mathcal{A}.$$

The second step is to prove that

$$\sup_{f \in \mathcal{A}} \sup_{\epsilon > 0} r_\epsilon^{-1} \mathbb{E}\|\hat{f} - f\|_2^2 < \infty.$$

From now on, we denote by $MS(\hat{f}, (r_\epsilon)_\epsilon)$ the maxiset of the procedure $\hat{f}$ associated with the rate of convergence $r = (r_\epsilon)_\epsilon$. The two steps to establish a maxiset result can be rewritten as the following embedding properties: $MS(\hat{f}, (r_\epsilon)_\epsilon) \subseteq \mathcal{A}$ corresponds to the first step and $\mathcal{A} \subseteq MS(\hat{f}, (r_\epsilon)_\epsilon)$ to the second one.

Although the maxiset approach is not extremely different from the minimax one, it is more optimistic since it provides a functional space directly connected to the estimation procedure. Thus this theoretical criterion to measure the performance of a chosen procedure appears to be more interesting for practical purposes. Indeed describing the maxiset of a procedure means knowing the entire functional space of well estimated functions. According to this point of view, the larger the maxiset, the better the procedure. Moreover it is interesting to remark that if a procedure $\hat{f}^*$ is $\mathcal{F}$-minimax optimal then

$$\mathcal{F} \subseteq MS(\hat{f}^*, (\rho_{\mathcal{F}}(\epsilon))_\epsilon).$$

In the wavelet setting and using the maxiset approach, many results have appeared in nonparametric statistics. Cohen, De Vore, Kerkyacharian and Picard [6] and Rivoirard [22] have proved that linear procedures are outperformed by non linear ones in the density estimation model and the white noise model. In particular, they have identified the maxisets of thresholding procedures with the intersection of Besov spaces and specific Lorentz spaces, called weak Besov spaces. More recently, Rivoirard [23] has shown that the maxisets of thresholding procedures coincide with those of classical Bayesian procedures associated with heavy tailed priors. Kerkyacharian and Picard [16] have proved that under some conditions, the maxiset of local bandwidth selection procedure is at least as large as the one of the hard thresholding procedure, but they have let two open questions: what is exactly the maxiset of the local bandwidth selection procedure? Does the local bandwidth selection procedure outperform the hard thresholding procedure in the maxiset sense?

The goal of this paper is to provide a new wavelet procedure which performs very well under both, the minimax and the maxiset approaches. In particular we aim at building a data-driven procedure which has better performances than the hard thresholding procedure. According to Autin [1] and [2] the only way to succeed in doing this is to consider procedures which are not *elitist*, i.e. that allow to use some empirical wavelet coefficients smaller than a threshold for the reconstruction of the signal. Here we propose a wavelet procedure (hard tree rule) inspired from the local bandwidth selection procedure of Lepski [17].

Firstly, this new wavelet procedure depends on the choice of a maximal scale $j_{max}$ to ensure the calculability of the estimate. According to this parameter, any empirical wavelet coefficient of the target function with a level index $j$ larger than or equal to $j_{max}$ will not be considered for the reconstruction. As in Autin [1] and [2], the choice of this maximal scale will have a direct consequence on the shape of the maxiset.

Secondly, the new procedure is based on thresholding methods associated with hereditary constraints (see Engel [13]). Using some ideas from tree approximation (see Cohen, Dahmen, Daubechies and De Vore [4], Engel [13]), from coding theory (see De Vore, Johnson, Pan and Sharpley [9], Said and Pearlman [24], Shapiro [25]), and from Image Processing (see Wainwright, Martin et al. [26] and Azimifar et al. [3]) we show that our new way of organizing the signal reconstruction allows to build a procedure with a very large maxiset. This new procedure outperforms the hard thresholding one as well as any elitist procedure in the maxiset sense.

The paper is organized as follows. Section 2 is devoted to the description of the model and the definitions of the basic tools we shall need. In Section 3, we describe the method to construct our wavelet procedure and we show the relationship with the local bandwidth selection procedure of Lepski. The minimax and maxiset performances of the procedure are studied in Sections 4 and 5. The comparison between the performances of this procedure and the hard thresholding ones is discussed in Section 6. A short conclusion is given in Section 7 while the proofs of our results are given in the Appendix.

## 2. Model and definitions

### 2.1. Model

We consider the white noise model: $X_\epsilon(.)$ is a random variable satisfying the following equation:

$$X_\epsilon(dt) = f(t)dt + \epsilon W(dt), \qquad t \in [0, 1[$$

where

- $0 < \epsilon < \frac{1}{2}$ is the noise level,
- $f$ is a function defined on $[0, 1[$,
- $W(.)$ is the standard Brownian on $[0, 1[$.

Let $\{\psi_{jk}(\cdot), j \geq -1, k \in \mathbb{Z}\}$ be a wavelet basis of $\mathbb{L}_2([0, 1[)$ with $N$ vanishing moments ($N \in \mathbb{N}^*$) built by multi resolution analysis from a scaling function and a wavelet supported on $[0, 2N - 1[$. Any $f \in \mathbb{L}_2([0, 1[)$ can be represented as:

$$f = \sum_{j \geq -1} \sum_{k \in \mathbb{Z}} \beta_{jk} \psi_{jk} = \sum_{j \geq -1} \sum_{k \in \mathbb{Z}} (f, \psi_{jk})_{L_2} \psi_{jk}. \tag{2.1}$$

There exists a constant $S_\psi$ such that at each level $j \geq -1$ there are less than or equal to $\mathcal{K}_\psi^{(j)} = 2^j \times S_\psi$ non-zero wavelet coefficients. Hence, at each level $j$, the sum over $k$ in (2.1) can be replaced by the sum over $k \in \mathcal{K}_\psi^{(j)}$.

In our setting, we can get all the observations: $y_{jk} = X_\epsilon(\psi_{jk}) = \beta_{jk} + \epsilon Z_{jk}$ where $Z_{jk}$ are independent Gaussian variables $\mathcal{N}(0, 1)$.

All along this paper, for a real $\eta \geq 1$, we write $2^{j_\lambda} \sim \lambda^{-2\eta}$ to design the integer $j_\lambda$ such that $2^{-j_\lambda} \leq \lambda^{2\eta} < 2^{1 - j_\lambda}$.

### *2.2. Definitions*

**Definition 2.1.** *We say that an interval $I_{jk}$ is dyadic if it corresponds to the support of the function $\psi_{jk}$ and we denote by $|I_{jk}| = l_\psi \, 2^{-j}$ its length (where $l_\psi$ is the size of the support of the mother wavelet function).*

**Definition 2.2.** *Let $\lambda > 0$ and $I_{jk}$ be a dyadic interval such that $0 \leq j < j_\lambda$. We denote by $\mathcal{T}_{jk}^{(\eta)}(\lambda)$ the binary tree containing the set of the dyadic intervals such that the following properties are satisfied:*

- *$I_{jk} \in \mathcal{T}_{jk}^{(\eta)}(\lambda)$.*
- *$I_{j'k'} \in \mathcal{T}_{jk}^{(\eta)}(\lambda) \Longrightarrow I_{j'k'} \subseteq I_{jk}$ and $|I_{j'k'}| > l_\psi \lambda^{2\eta}$.*
- *Two distinct dyadic intervals of $\mathcal{T}_{jk}^{(\eta)}(\lambda)$ with same length have their interiors disjointed.*
- *The numbers of dyadic intervals of $\mathcal{T}_{jk}^{(\eta)}(\lambda)$ of length $l_\psi 2^{-j'}$ ($j \leq j' < j_\lambda$) is equal to $2^{j'-j}$.*
- *Any set of all dyadic intervals of $\mathcal{T}_{jk}^{(\eta)}(\lambda)$ with same length is forming a partition of $I_{jk}$.*

## 3. Construction of a new adaptive procedure

The aim of this section is to provide a new wavelet procedure based on thresholding methods which takes advantage on the dyadic structure of the wavelet decomposition.

Let

$$f(t) = \sum_{j=-1}^{\infty} \sum_{k \in \mathcal{K}_\psi^{(j)}} \beta_{jk} \psi_{jk}(t), \quad t \in [0, 1[$$

be a function to be estimated from the observations $y_{jk}$ of its wavelet coefficients $\beta_{jk}$. We propose to estimate the function only using a finite number of observations of wavelet coefficients, that's why we consider the following family of Keep-Or-Kill estimators:

$$\mathcal{F}_K(\epsilon) = \left\{ \tilde{f}(.) = \sum_{j=-1}^{j_{max}(\epsilon)-1} \sum_{k \in \mathcal{K}_\psi^{(j)}} \gamma_{jk} y_{jk} \psi_{jk}(.), \ \gamma_{jk} \in \{0, 1\} \right\}.$$

Any procedure in $\mathcal{F}_K(\epsilon)$ does not use the empirical wavelet coefficients $y_{jk}$ for which the level $j$ is larger than or equal to $j_{max}(\epsilon)$. This condition ensures that any procedure of $\mathcal{F}_K(\epsilon)$ is numerically calculable. As we shall see in Section 5, on the choice of the maximum scale $j_{max}$ will depend the maxiset of the procedure considered.

In the sequel, we shall set $\lambda_\epsilon = m\epsilon\sqrt{\log(\epsilon^{-1})}$ where $m$ is an absolute constant which will be chosen later and, for a fixed real number $\eta \geq 1$ (maximum scale parameter), we shall denote by $j_{\lambda_\epsilon}$ the integer such that $2^{j_{\lambda_\epsilon}} \sim \lambda_\epsilon^{-2\eta}$ and we shall put $j_{max}(\epsilon) = j_{\lambda_\epsilon}$.

### 3.1. Definition of the hard tree procedure

Let us consider the following procedure, namely the hard tree procedure, defined for $\eta \geq 1$ by:

$$\tilde{f}_T(.) = \sum_{k \in \mathcal{K}_\psi^{(-1)}} y_{-1k} \psi_{-1k}(.) + \sum_{j=0}^{j_{\lambda_\epsilon}-1} \sum_{k \in \mathcal{K}_\psi^{(j)}} \gamma_{jk} y_{jk} \psi_{jk}(.) \tag{3.1}$$

with

- $\gamma_{jk} = 1$ if there exists $I_{j'k'}$ in $\mathcal{T}_{jk}^{(\eta)}(\lambda_\epsilon)$ such that $|y_{j'k'}| > \lambda_\epsilon$,
- $\gamma_{jk} = 0$ otherwise.

At first glance, this estimator is not very different from the hard thresholding one recalled in (4.2). It consists in keeping the empirical coefficients larger than $\lambda_\epsilon$ and somehow, "in filling the holes" in each binary tree $\mathcal{T}_{jk}(\lambda_\epsilon)$, as we can see in Figure 1.

Notice that the hard tree estimator minimizes a penalized criterion. Indeed,

$$\tilde{f}_T = Arg \min_{\hat{f} \in \mathcal{F}_K(\epsilon)} \sum_{j=0}^{j_{\lambda_\epsilon}-1} \sum_{k \in \mathcal{K}_\psi^{(j)}} (\gamma_{jk} - 1)^2 |\bar{y}_{jk}(\lambda_\epsilon)|^2 + \gamma_{jk}^2 \lambda_\epsilon^2$$

where $|\bar{y}_{jk}(\lambda_\epsilon)| := \max\{|y_{j'k'}|, \ I_{j'k'} \in \mathcal{T}_{jk}^{(\eta)}(\lambda_\epsilon)\}$.

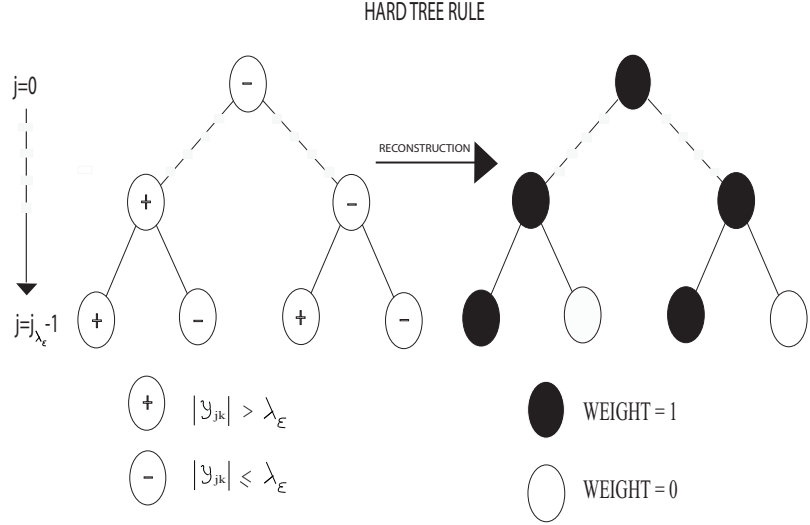

Fɪɢ 1.

Moreover, this procedure is a tree rule (Engel [13]) since it satisfies the following hereditary constraints:

$$\gamma_{jk} = 1 \implies \forall I_{j'k'} \ / \ I_{jk} \in \mathcal{T}_{j'k'}^{(\eta)}(\lambda_\epsilon), \quad \gamma_{j'k'} = 1,$$
$$\gamma_{jk} = 0 \implies \forall I_{j'k'} \in \mathcal{T}_{jk}^{(\eta)}(\lambda_\epsilon), \quad \gamma_{j'k'} = 0.$$

Tree-structures are often used in approximation theory and coding theory. For more details, we refer the reader to the papers of Cohen, Dahmen, Daubechies and De Vore [4], Cohen, Daubechies, Guleryuz and Orchard [5], De Vore, Johnson, Pan and Sharpley [9], Said and Pearlman [24] and Shapiro [25].

### 3.2. Algorithm for the construction of hard tree rule

In this paragraph, we give the method to construct the hard tree procedure, assuming that the noise level $\epsilon$ is known.

*Algorithm*
*Setup*:

- Choose the reals $\eta \geq 1$ and $m > 0$ and put $\lambda_\epsilon = m\epsilon\sqrt{\log(\epsilon^{-1})}$;
- Identify $j_{\lambda_\epsilon} = min\{j \in \mathbb{N}, 2^j \geq \lambda_\epsilon^{-2\eta}\}$.

*Construction steps*:

- Compute $y_{jk}$ with $k \in \mathcal{K}_\psi^{(j)}$ and level $j < j_{\lambda_\epsilon}$;

- Threshold any $y_{jk}$ at level $\lambda_\epsilon$ and construct the set of indices

$$\mathcal{I}(\lambda_\epsilon) = \{(j,k),\ j \geq 0,\ k \in \mathcal{K}_\psi^{(j)},\ |y_{jk}| > \lambda_\epsilon\};$$

- Construct the set of indices

$$\mathcal{S}(\lambda_\epsilon) = \bigcup_{(j,k) \in \mathcal{I}(\lambda_\epsilon)} \{(j',k'),\ \text{with } I_{j'k'} \supseteq I_{jk} \text{ and } I_{jk} \in \mathcal{T}_{j'k'}^{(\eta)}(\lambda_\epsilon)\}.$$

*Return*:

- The estimator $\tilde{f}_T = \sum\limits_{k \in \mathcal{K}_\psi^{(-1)}} y_{-1k} \psi_{-1k} + \sum\limits_{(j,k) \in \mathcal{S}(\lambda_\epsilon)} y_{jk} \psi_{jk}$.

### 3.3. Connection with Lepski's rule

In this paragraph, we show that the hard tree rule can be viewed as a *wavelet-version* of the bandwidth selection procedure of Lepski [17] when the chosen wavelet basis is the Haar one.

Notice that, in the Haar case, any dyadic interval is on the form $I_{jk} = [\frac{k}{2^j}, \frac{k+1}{2^j}[$ with $j \in \mathbb{N}$ and $k \in \{0,\dots,2^j-1\}$. Moreover, one gets a characterization of its wavelet components $\psi_{jk}(.)$. Indeed

$$\psi_{-1k}(.) = \psi_{-1}(.-k) \quad \text{and} \quad \psi_{jk}(.) = 2^{\frac{j}{2}}\psi(2^j.-k),\ \text{with}$$
$$\psi_{-1}(.) = \mathbf{1}_{[0,1[}(.) \quad \text{and} \quad \psi(.) = \mathbf{1}_{[0,\frac{1}{2}[}(.) - \mathbf{1}_{[\frac{1}{2},1[}(.).$$

With this particular choice of wavelet basis, the hard tree procedure is defined by

$$\tilde{f}_T(.) = y_{-10} + \sum_{j=0}^{j_{\lambda_\epsilon}-1} \sum_{k \in \mathcal{K}_\psi^{(j)}} \gamma_{jk} y_{jk} \left( \mathbf{1}_{[\frac{k}{2^j},\frac{2k+1}{2^{j+1}}[}(.) - \mathbf{1}_{[\frac{2k+1}{2^{j+1}},\frac{k+1}{2^j}[}(.) \right),$$

with $\gamma_{jk} = 2^{\frac{j}{2}}$ if there exists $I_{j'k'} \subseteq I_{jk}$ such that $|I_{j'k'}| > \lambda_\epsilon^{2\eta}$ and $|y_{j'k'}| > \lambda_\epsilon$; $\gamma_{jk} = 0$ otherwise.

Let us now briefly recall the definition of the local bandwidth selection rule (see Lepski [17] or Lepski, Mammen and Spokoiny [18] for more details).

*Local bandwidth selection rule*

Let $K$ be a compactly supported bounded kernel such that $\|K\|_{\mathbb{L}_2} = 1$. For any $j \in \mathbb{N}$ and any $(t,u) \in [0,1]^2$, let us denote

$$K_j(t,u) = 2^j K(2^j t, 2^j u) \text{ and } \hat{K}_j(t) = \int_0^1 K_j(t,u) dX_\epsilon(u).$$

Let us define the index $\hat{j}(t)$ as the minimum of *admissible* $j$'s at the point $t$, where $j < j_{\lambda_\epsilon}$ is admissible at the point $t$ if $j = j_{\lambda_\epsilon}$ or

$$|\hat{K}_{j'+1}(t) - \hat{K}_{j'}(t)| \leq 2^{\frac{j'}{2}} \lambda_\epsilon \qquad \forall j \leq j' < j_{\lambda_\epsilon}. \tag{3.2}$$

The local bandwidth selection estimator $\hat{f}_L$ is defined by:

$$\hat{f}_L(t) = \hat{K}_{\hat{j}(t)}(t).$$

The definition of the hard tree rule is close to the definition of the local bandwidth selection procedure. Indeed, let us adapt the notion of *admissibility* from kernel estimators to wavelet estimators by considering the family of estimators $(\hat{f}_j)_{j \in \mathbb{N}}$ defined as follows:

- $\hat{f}_0(t) = y_{-10}\psi_{-10}(t)$

- $\hat{f}_{j+1}(t) = \hat{f}_j(t) + \sum_{k=0}^{2^j-1} y_{jk}\psi_{jk}(t).$

If for any $t \in [0, 1[$ we denote by $I_j^t$ the dyadic interval containing $t$ such that $|I_j^t| = 2^{-j}$, then

$$|\hat{f}_{j+1}(t) - \hat{f}_j(t)| = \left| \sum_{k=0}^{2^j-1} y_{jk}\psi_{jk}(t) \right| = 2^{\frac{j}{2}}|y_{I_j^t}|. \tag{3.3}$$

**Definition 3.1.** *We say that an integer $j$ is $(t,T)$-admissible if:*
*either $j = j_{\lambda_\epsilon}$ or, for all $j \le j' < j_{\lambda_\epsilon}$, for all $t' \in I_j^t$:*

$$|\hat{f}_{j'+1}(t') - \hat{f}_{j'}(t')| \le 2^{\frac{j'}{2}}\lambda_\epsilon.$$

Denote $\hat{\hat{j}}_T(t) = \inf\{j;\ j \text{ is (t,T)-admissible}\}$. Still using (3.3) we can observe that:

$$\hat{f}_{\hat{j}_T(t)}(t) = \tilde{f}_T(t). \tag{3.4}$$

So, by adapting the notion of *admissibility* from kernel procedures to wavelet procedures, we have shown that the adaptive procedure (hard tree rule) and Lepski's rule are analogous when considering the particular choice of the kernel $K$:

$$K_j(x, y) = 2^j \sum_{k=0}^{2^j-1} \psi_{-1}(2^j x - k)\psi_{-1}(2^j y - k).$$

## 4. Minimax result

In this paragraph we aim at studying the performance associated with the hard tree rule in the minimax context.

At first, let us recall the definition of Besov spaces $\mathcal{B}_{2,\infty}^s$, with $0 < s < N$.

**Definition 4.1.** *Let $0 < s < N$. We say that a function $f \in \mathbb{L}_2([0, 1[)$ belongs to the Besov space $\mathcal{B}_{2,\infty}^s$ if and only if:*

$$\sup_{J \ge 0} 2^{2Js} \sum_{j=J}^{\infty} \sum_{k \in \mathcal{K}_\psi^{(j)}} \beta_{jk}^2\ < \infty.$$

Besov spaces constitute a large class of functional spaces. Recall that for any $s > 0$ Sobolev space $H^s$ is included in $\mathcal{B}_{2,\infty}^s$. Moreover, if $\subsetneq$ denotes the strict inclusion between two functional spaces,

$$\mathcal{B}_{2,\infty}^s \subsetneq \mathcal{B}_{2,2}^{s'} \quad \text{for any } 0 < s' < s < N. \tag{4.1}$$

Besov spaces are important in statistics since the maximal spaces of many classical procedures like linear procedures (see Kerkyacharian and Picard [14] and Rivoirard [22]) and thresholding procedures (see Cohen, De Vore, Kerkyacharian and Picard [6] and Kerkyacharian and Picard [16]) are included in Besov spaces.

We prove in the following theorem that the hard tree procedure is $\mathcal{B}_{2,\infty}^s$-minimax optimal up to a logarithmic term which is known to be the price to pay for adaptation.

**Theorem 4.1.** *Let $0 < s < N$ and $\eta \geq 1$. Choose $m \geq 4\sqrt{3\eta}$. Then for any $f \in \mathcal{B}_{2,\infty}^s$*

$$\sup_{\epsilon > 0} \lambda_\epsilon^{-4s/(1+2s)} \mathbb{E} \|\tilde{f}_T - f\|_2^2 < \infty.$$

This result is just a consequence of Theorem 5.2 using the embedding properties (5.1) and (5.2). This theorem shows that the hard tree procedure described in Section 3 performs very well. Moreover, let us recall the minimax result for the hard thresholding procedure:

$$\tilde{f}_H(.) = \sum_{k \in \mathcal{K}_\psi^{(-1)}} y_{-1k}\psi_{-1k}(.) + \sum_{j=0}^{j_{\lambda_\epsilon}-1} \sum_{k \in \mathcal{K}_\psi^{(j)}} y_{jk}\mathbf{1}\{|y_{jk}| > \lambda_\epsilon\}\psi_{jk}(.). \tag{4.2}$$

**Theorem 4.2.** *Let $0 < s < N$ and $\eta \geq 1$. Choose $m \geq 4\sqrt{2\eta}$. Then for any $f \in \mathcal{B}_{2,\infty}^s$*

$$\sup_{\epsilon > 0} \lambda_\epsilon^{-4s/(1+2s)} \mathbb{E} \|\tilde{f}_H - f\|_2^2 < \infty.$$

This minimax result is a direct consequence of Theorem 5.1 of Section 5 using the embedding property (5.1).

**Remark 4.1.** *It is important to notice here that the minimax results given in Theorems 4.1 and 4.2 are valid for any choice of compactly supported wavelet basis provided that its number of vanishing moments $N$ is strictly greater than $s$.*

Following the two last theorems

**Corollary 4.1.** *For any $0 < s < N$ and any choice of $\eta \geq 1$, the hard tree procedure has the same performance as the hard thresholding procedure from the minimax point of view when considering the same threshold level $\lambda_\epsilon = m\epsilon\sqrt{\log(\epsilon^{-1})}$ with $m \geq 4\sqrt{3\eta}$. Precisely, both procedures are $\mathcal{B}_{2,\infty}^s$-minimax optimal (up to a logarithmic term).*

A natural question arises here: could these procedures be discriminated when adopting the maxiset point of view? The answer is *YES* as we shall see.

## 5. Maxiset result

In this section, we aim at calculating the maxiset associated with the hard tree procedure so as to compare it with the one of the hard thresholding procedure when the rate of convergence is $(\lambda_\epsilon^{4s/(1+2s)})_\epsilon$, $0 < s < N$.

At first we propose to recall the maxiset result given by Kerkyacharian and Picard [15] for the hard thresholding estimator.

### 5.1. Maxiset of the hard thresholding procedure

Let us introduce the following functional space.

**Definition 5.1.** *Let $0 < r < 2$. We say that a function $f$ belongs to the weak Besov space $W_r$ if and only if:*

$$\sup_{\lambda > 0} \lambda^{r-2} \sum_{j=0}^{\infty} \sum_{k \in \mathcal{K}_\varphi^{(j)}} \beta_{jk}^2 \mathbf{1}\{|\beta_{jk}| \leq \lambda\} \ < \infty.$$

Weak Besov spaces compose a sub-family of Lorentz spaces (see Lorentz [19], [20] or De Vore and Lorentz [11]). There exists a natural relationship between Besov spaces and weak Besov spaces. The following embedding can be easily proved (see for instance Rivoirard [22]):

$$\mathcal{B}_{2,\infty}^s \subsetneq \mathcal{B}_{2,\infty}^{\frac{s}{\eta(1+2s)}} \cap W_{\frac{2}{1+2s}} \quad \text{for any } 0 < s < N \text{ and any } \eta \geq 1. \tag{5.1}$$

Kerkyacharian and Picard [15] and [16] have pointed out the strong connection between these functional spaces and the hard thresholding procedure.

**Theorem 5.1** (Kerkyacharian-Picard). *Let $0 < s < N$ and $\eta \geq 1$. For any $m \geq 4\sqrt{2\eta}$, we have the following equivalence:*

$$\sup_{0 < \epsilon < \frac{1}{2}} \lambda_\epsilon^{-\frac{4s}{1+2s}} \ \mathbb{E}\|\tilde{f}_H - f\|_2^2 < \infty \iff f \in \mathcal{B}_{2,\infty}^{\frac{s}{\eta(1+2s)}} \cap W_{\frac{2}{1+2s}},$$

*that is to say, using the maxiset notation:*

$$MS(\tilde{f}_H, (\lambda_\epsilon^{\frac{4s}{1+2s}})_\epsilon) = \mathcal{B}_{2,\infty}^{\frac{s}{\eta(1+2s)}} \cap W_{\frac{2}{1+2s}}.$$

### 5.2. Maxiset of the hard tree procedure

In this paragraph, we exhibit the maxiset associated with the hard tree procedure associated with the rate $(\lambda_\epsilon^{4s/(1+2s)})_\epsilon$.

Let us first define another functional space that will be useful in the characterization of the maximal space associated with hard tree procedure.

**Definition 5.2.** *Let $0 < r < 2$ and $\eta \geq 1$. We say that a function $f$ belongs to the space $W_{r,\eta}^T$ if and only if:*

$$\sup_{\lambda > 0} \lambda^{r-2} \sum_{j=0}^{j_\lambda - 1} \sum_{k \in \mathcal{K}_\psi^{(j)}} \beta_{jk}^2 \mathbf{1} \left\{ \forall I_{j'k'} \in \mathcal{T}_{jk}^{(\eta)}(\lambda), |\beta_{j'k'}| \leq \frac{\lambda}{2} \right\} < \infty.$$

In contrast of weak Besov spaces, note that the spaces $W_{r,\eta}^T$ ($0 < r < 2, \eta \geq 1$) are not invariant under permutations of wavelet coefficients within each scale.

The following proposition shows that, for the same parameter $r$ ($0 < r < 2$), any functional space $W_{r,\eta}^T$ contains the weak Besov space $W_r$. Thanks to this result, a comparison between the maximal sets of hard tree rule and the hard thresholding rule will be possible, as we will see in Section 6.

**Proposition 5.1.** *For any $0 < r < 2$ and any $\eta \geq 1$, we have the following inclusion spaces:*

$$W_r \quad \subsetneq \quad W_{r,\eta}^T. \tag{5.2}$$

Proposition 5.1 shows that for any parameters $0 < r < 2$ and $\eta \geq 1$, spaces $W_r$ and $W_{r,\eta}^T$ are *different*.

**Theorem 5.2.** *Let $0 < s < N$ and $\eta \geq 1$. For any $m \geq 4\sqrt{3\eta}$, we have the following equivalence:*

$$\sup_{0 < \epsilon < \frac{1}{2}} \lambda_\epsilon^{-\frac{4s}{1+2s}} \mathbb{E} \|\tilde{f}_T - f\|_2^2 < \infty \iff f \in \mathcal{B}_{2,\infty}^{\overline{\frac{s}{\eta(1+2s)}}} \cap W_{\frac{2}{1+2s},\eta}^T,$$

*that is to say, using the maxiset notation:*

$$MS(\tilde{f}_T, (\lambda_\epsilon^{\frac{4s}{1+2s}})_\epsilon) = \mathcal{B}_{2,\infty}^{\overline{\frac{s}{\eta(1+2s)}}} \cap W_{\frac{2}{1+2s},\eta}^T.$$

To prove this theorem we shall need the following proposition.

**Proposition 5.2.** *Fix $0 < \lambda_0 < 1$ and $\eta \geq 1$. For any $0 < r < 2$ and any $f \in \mathcal{B}_{2,\infty}^{(2-r)/4\eta} \cap W_{r,\eta}^T$, then*

$$\sup_{0 < \lambda \leq \lambda_0} \lambda^r \left[ \log\left(\frac{1}{\lambda}\right) \right]^{-1} \sum_{j=0}^{j_\lambda - 1} \sum_{k \in \mathcal{K}_\psi^{(j)}} \mathbf{1} \left\{ \exists I_{j'k'} \in \mathcal{T}_{jk}^{(\eta)}(\lambda) \, / \, |\beta_{j'k'}| > \frac{\lambda}{2} \right\} < \infty. \tag{5.3}$$

**Remark 5.1.** *It is important to notice here that the maxiset results given in Theorems 5.1 and 5.2 are valid for any choice of compactly supported wavelet basis provided that its number of vanishing moments $N$ is strictly greater than $s$.*

Following the two previous sections, let us comment the minimax and maxiset performances of the hard tree rule.

## 6. On the performances of the hard tree procedure

### 6.1. Consequences of previous results

Judging from Corollary 4.1 of Section 4, the hard tree procedure and the hard thresholding one are equivalent in the minimax sense.

According to Proposition 5.1 and Theorem 5.2 we easily deduce that the hard tree procedure performs very well in the maxiset sense. Indeed, for a chosen $\eta \geq 1$, its maxiset for the rate $(\lambda_\epsilon^{4s/(1+2s)})_\epsilon$ corresponds to the intersection between the usual Besov space $\mathcal{B}_{2,\infty}^{\overline{s/(1+2s)}}$ and another functional space $W_{\frac{2}{1+2s},\eta}^T$ strictly larger than the classical weak Besov space $W_{\frac{2}{1+2s}}$. Hence,

**Corollary 6.1.** *In the maxiset sense, the hard tree procedure is at least as good as the hard thresholding procedure since its maxiset for the rate $(\lambda_\epsilon^{4s/(1+2s)})_\epsilon$ contains the hard thresholding procedure one.*

It is important to notice that a strict inclusion between the maxisets of the hard tree rule and the hard thresholding rule can not be immediately deduced from previous results because of the intersections with the Besov space. At present, it is an *open question* whether the inclusion between maxisets is strict or not. Nevertheless we give in the sequel results which address a slightly weaker problem.

### 6.2. More results on spaces embeddings

**Proposition 6.1.** *For any $0 < s < N$ and any $\eta \geq 1$ the following spaces embedding holds:*

$$\mathcal{B}_{2,\infty}^{\overline{\frac{s}{\eta(1+2s)}},*} \cap W_{\frac{2}{1+2s}} \quad \subsetneq \quad \mathcal{B}_{2,\infty}^{\overline{\frac{s}{\eta(1+2s)}},*} \cap W_{\frac{2}{1+2s},\eta}^T,$$

*with*

$$\mathcal{B}_{2,\infty}^{u,*} = \left\{ f = \sum_{j=-1}^{\infty} \sum_{k \in \mathcal{K}_\psi^{(j)}} \beta_{jk} \psi_{jk}; \ \sup_{J \geq -1} J^{-1} . 2^{2Ju} \sum_{j=J}^{\infty} \sum_{k \in \mathcal{K}_\psi^{(j)}}^{2^j - 1} \beta_{jk}^2 \ < \infty \right\}.$$

According to Proposition 6.1, the strict inclusions of functional spaces are still valid when intersecting $W_{\frac{2}{1+2s}}$ and $W_{\frac{2}{1+2s},\eta}^T$ with the hybrid Besov space $\mathcal{B}_{2,\infty}^{\overline{\frac{s}{1+2s}},*}$. From this result one immediately derives

**Corollary 6.2.** *Let $0 < s < N$ and $\eta \geq 1$. The following spaces embedding holds for all $u < \frac{s}{\eta(1+2s)}$:*

$$\mathcal{B}_{2,\infty}^u \cap W_{\frac{2}{1+2s}} \quad \subsetneq \mathcal{B}_{2,\infty}^u \cap W_{\frac{2}{1+2s},\eta}^T.$$

*Moreover,*

$$\bigcap_{u<\frac{s}{\eta(1+2s)}} \mathcal{B}^u_{2,\infty} \cap W_{\frac{2}{1+2s}} \quad \subsetneq \bigcap_{u<\frac{s}{\eta(1+2s)}} \mathcal{B}^u_{2,\infty} \cap W^T_{\frac{2}{1+2s},\eta}.$$

Also, the embeddings of spaces with strict inclusion are still valid when considering intersection of spaces very close to the maxisets we have studied. Hence it is reasonable to claim that hard tree procedure is better than the hard thresholding procedure in the maxiset sense.

### *6.3. On the choice of parameter $\eta$*

In Sections 4 and 5 we gave minimax and maxiset results on the hard tree procedure for any choice of parameter $\eta \geq 1$. Precisely, the regularity parameter of the Besov space appearing in the maxisets of the hard tree and hard thresholding procedures depends on the choice of $\eta$. Hence it could be interested to know if an optimal choice of $\eta$ could be possible so as to build the hard tree rule with the largest maxiset. In fact, there is no doubt that the bigger the parameter $\eta$ the larger the maxiset of the hard tree rule. Indeed

**Proposition 6.2.** *For any $0 < s < N$ and any $1 \leq \eta_1 < \eta_2$, the following spaces embeddings hold:*

$$\mathcal{B}^{\overline{\frac{s}{\eta_1(1+2s)}}}_{2,\infty} \cap W^T_{\frac{2}{1+2s},\eta_1} \subsetneq \mathcal{B}^{\overline{\frac{s}{\eta_2(1+2s)}}}_{2,\infty} \cap W^T_{\frac{2}{1+2s},\eta_2}.$$

Nevertheless our results are asymptotic. In fact, if at first glance we opt for a choice of a very large $\eta$, we must be careful to the change for the worse of rate of convergence considered. Indeed, choosing a large $\eta$ implies taking a large $m$. As a consequence the rate of convergence $(\lambda^{4s/(1+2s)}_\epsilon)_\epsilon$ goes more slowly for such a choice.

## 7. Conclusion

The key point of this paper was to prove that a way to build very performing procedures is to combine thresholding methods and tree structure. Indeed the maxiset of the new wavelet procedure called hard tree procedure is proved to perform very well in the minimax and the maxiset settings. Although this procedure looks like the hybrid version of Lepski's procedure proposed by Picard and Tribouley [21], namely hard stem rule, this one is different (see Autin [1]) and presents more advantages comparing to the hard stem rule. Firstly, Autin [1] has proved that the maxiset of the hard tree rule contains the one of the hard stem rule. Secondly the hard stem rule is a procedure which is only defined with the Haar wavelet basis. Indeed, the hard stem procedure is built at fixed $t \in [0,1[$ and therefore especially requires wavelet functions $\psi_{jk}$ with disjoint supports. Here the hard tree rule is defined for any compactly supported wavelet basis.

## 8. Appendix

### 8.1. Proofs of Propositions

Here and later, the constants $C$ represent all the constants we shall need and can be different from one line to one other.

*Proof of Proposition 5.1.* Let $\eta \geq 1$. The large inclusion is obvious when remarking that for any sequence of wavelet coefficients $(\beta_{jk}, j \geq 0, k)$, for any $0 < \lambda < 1$ and any $0 \leq j < j_\lambda$

$$\sum_k \beta_{jk}^2 \mathbf{1}\left\{\forall I_{j'k'} \in \mathcal{T}_{jk}^{(\eta)}(\lambda), |\beta_{j'k'}| \leq \frac{\lambda}{2}\right\} \leq \sum_k \beta_{jk}^2 \mathbf{1}\left\{|\beta_{jk}| \leq \frac{\lambda}{2}\right\}.$$

The strict inclusion is a direct consequence of Proposition 6.1. □

*Proof of Proposition 5.2.* Let $f \in \mathcal{B}_{2,\infty}^{(2-r)/4\eta} \cap W_{r,\eta}^T$ and $0 < \lambda \leq \lambda_0$. We set for any $u \in \mathbb{N}$, $2^{j_{\lambda,u}} \sim (2^{1+u}\lambda)^{-2\eta}$. We have

$$\sum_{j=0}^{j_\lambda-1} \sum_k \mathbf{1}\left\{\exists I_{j'k'} \in \mathcal{T}_{jk}^{(\eta)}(\lambda) \,/\, |\beta_{j'k'}| > \frac{\lambda}{2}\right\}$$

$$\leq \quad \sum_{u \geq 0} \sum_{j=0}^{j_\lambda-1} \sum_k (j+1)\mathbf{1}\{2^{u-1}\lambda < |\beta_{jk}| \leq 2^u\lambda, \ \ \forall I_{j'k'} \in \mathcal{T}_{jk}^{(\eta)}(\lambda), |\beta_{j'k'}| \leq 2^u\lambda\}$$

$$\leq \quad C \, j_\lambda \sum_{u \geq 0} (2^{u-1}\lambda)^{-2}$$

$$\times \sum_{j=0}^{j_\lambda-1} \sum_k \beta_{jk}^2 \mathbf{1}\{2^{u-1}\lambda < |\beta_{jk}| \leq 2^u\lambda, \ \forall I_{j'k'} \in \mathcal{T}_{jk}^{(\eta)}(2^{1+u}\lambda), |\beta_{j'k'}| \leq 2^u\lambda\}$$

$$\leq \quad C \, \log\left(\frac{1}{\lambda}\right) \sum_{u \geq 0} (2^{u-1}\lambda)^{-2}$$

$$\times \sum_{j=0}^{j_{\lambda,u}-1} \sum_k \beta_{jk}^2 \mathbf{1}\{2^{u-1}\lambda < |\beta_{jk}| \leq 2^u\lambda, \ \forall I_{j'k'} \in \mathcal{T}_{jk}^{(\eta)}(2^{1+u}\lambda), |\beta_{j'k'}| \leq 2^u\lambda\}$$

$$+ \, C \log\left(\frac{1}{\lambda}\right) \sum_{u \geq 0} (2^{u-1}\lambda)^{-2} \sum_{j=j_{\lambda,u}}^{\infty} \sum_k \beta_{jk}^2$$

$$\leq \quad C \, \log\left(\frac{1}{\lambda}\right) \lambda^{-r}.$$

The last inequality uses the fact that $f \in \mathcal{B}_{2,\infty}^{(2-r)/4\eta} \cap W_{r,\eta}^T$. □

To prove some strict inclusions in Proposition 6.1 and Proposition 6.2, let us introduce the function $h[m, \alpha, \alpha_1, \alpha_2](.)$ of $\mathbb{L}_2([0, 1[)$ for which the sequence of Haar wavelet coefficients $(\beta_{jk})_{jk}$ satisfies at each level j:

- if $j$ is even, $\lfloor (mj+1)2^{j\alpha} \rfloor \wedge 2^j$ wavelet coefficients $\beta_{jk}$ are equal to $2^{-\alpha_1 j}$, the others are equal to 0,
- if $j$ is odd, $\lfloor (mj+1)2^{j\alpha} \rfloor \wedge 2^j$ wavelet coefficients $\beta_{jk}$ are equal to $2^{-\alpha_2 j}$, the others are equal to 0,

and

$$\beta_{jk} \neq 0 \implies max(\beta_{j+1\ 2k}, \beta_{j+1\ 2k+1}) \neq 0.$$

*Proof of Proposition 6.1.* Fix $\eta > 1$ and $0 < s < N$. Looking at Proposition 5.1 with $r = 2(1+2s)^{-1}$, the large inclusion is obvious. For any $0 < s < 1$, the strict inclusion is given by considering the function $h[m, \alpha, \alpha_1, \alpha_2](.)$ with the parameters

$$m = 1,\ \alpha = (\eta(1+2s))^{-1},\ \alpha_1 = 1,\ \alpha_2 = (2\eta)^{-1},$$

which belongs to the space $\mathcal{B}_{2,\infty}^{\overline{\frac{s}{\eta(1+2s)}},*} \cap W_{\frac{2}{1+2s},\eta}^T$ but does not belong to the space $W_{\frac{2}{1+2s}}$. $\qquad\square$

*Proof of Proposition 6.2.* Let us first prove that for any $1 \leq \eta_1 < \eta_2$ and any $0 < s < N$,

$$\mathcal{B}_{2,\infty}^{\overline{\eta_1(1+2s)}} \cap W_{\frac{2}{1+2s},\eta_1}^T \subseteq W_{\frac{2}{1+2s},\eta_2}^T.$$

Let $f \in \mathcal{B}_{2,\infty}^{\overline{\eta_1(1+2s)}} \cap W_{\frac{2}{1+2s},\eta_1}^T$ and $0 < \lambda < 1$. We set $2^{j_{\lambda,1}} \sim \lambda^{-2\eta_1}$ and $2^{j_{\lambda,2}} \sim \lambda^{-2\eta_2}$. For any $j \geq 0$ and any $k$, we have

$$\sum_{j=0}^{j_{\lambda,2}-1} \sum_k \beta_{jk}^2 \mathbf{1}\left\{ \forall I_{j'k'} \in \mathcal{T}_{jk}^{(\eta_2)}(\lambda),\ |\beta_{j'k'}| \leq \frac{\lambda}{2} \right\}$$

$$\leq \sum_{j=0}^{j_{\lambda,1}-1} \sum_k \beta_{jk}^2 \mathbf{1}\left\{ \forall I_{j'k'} \in \mathcal{T}_{jk}^{(\eta_2)}(\lambda),\ |\beta_{j'k'}| \leq \frac{\lambda}{2} \right\} + \sum_{j=j_{\lambda,1}}^{\infty} \sum_k \beta_{jk}^2$$

$$\leq \sum_{j=0}^{j_{\lambda,1}-1} \sum_k \beta_{jk}^2 \mathbf{1}\left\{ \forall I_{j'k'} \in \mathcal{T}_{jk}^{(\eta_1)}(\lambda),\ |\beta_{j'k'}| \leq \frac{\lambda}{2} \right\} + \sum_{j=j_{\lambda,1}}^{\infty} \sum_k \beta_{jk}^2$$

$$\leq C\left( \lambda^{\frac{4s}{1+2s}} + 2^{-\frac{2sj_{\lambda,1}}{\eta_1(1+2s)}} \right)$$

$$\leq C\,\lambda^{\frac{4s}{1+2s}}.$$

The last inequality uses the fact that $f \in \mathcal{B}_{2,\infty}^{\overline{\frac{s}{\eta_1(1+2s)}}} \cap W_{\frac{2}{1+2s},\eta_1}^T$. Hence $f \in W_{\frac{2}{1+2s},\eta_2}^T$.

For any $0 < s < 1$, the strict inclusion is given by considering the function $h[m, \alpha, \alpha_1, \alpha_2](.)$ with the parameters

$$m = 0,\ \alpha = (\eta_2(1+2s))^{-1},\ \alpha_1 = \alpha_2 = (2\eta_2)^{-1},$$

which belongs to the space $\mathcal{B}_{2,\infty}^{\overline{\frac{s}{\eta_2(1+2s)}}} \cap W_{\frac{2}{1+2s},\eta_2}^T$ but does not belong to the space $\mathcal{B}_{2,\infty}^{\overline{\frac{s}{\eta_1(1+2s)}}}$. $\qquad\square$

### 8.2. *Proof of Theorem 5.2*

It needs two steps. At first we have to prove

$$STEP \ 1: \quad MS(\tilde{f}_T, (\lambda_\epsilon^{\frac{4s}{1+2s}})_\epsilon) \subseteq \mathcal{B}_{2,\infty}^{\overline{\eta(1+2s)}} \cap W_{\frac{2}{1+2s},\eta}^T. \tag{8.1}$$

Then

$$STEP \ 2: \quad MS(\tilde{f}_T, (\lambda_\epsilon^{\frac{4s}{1+2s}})_\epsilon) \supseteq \mathcal{B}_{2,\infty}^{\overline{\eta(1+2s)}} \cap W_{\frac{2}{1+2s},\eta}^T \tag{8.2}$$

must be proved.

*STEP 1:* Let $f \in MS(\tilde{f}_T, (\lambda_\epsilon^{\frac{4s}{1+2s}})_\epsilon)$. We have,

$$\sum_{j=j_{\lambda_\epsilon}}^{\infty} \sum_k \beta_{jk}^2 \leq \mathbb{E}\|\tilde{f}_T - f\|_2^2 \leq C\lambda_\epsilon^{\frac{4s}{1+2s}} \leq C 2^{-\frac{2j_{\lambda_\epsilon}s}{\eta(1+2s)}}.$$

So, using the continuity of $\lambda_\epsilon$ in 0, we deduce that

$$\sup_{J \geq -1} 2^{\frac{2Js}{\eta(1+2s)}} \sum_{j=J}^{\infty} \sum_k \beta_{jk}^2 \ < \infty.$$

It comes that $f \in \mathcal{B}_{2,\infty}^{\overline{\eta(1+2s)}}$. Let us now denote for any $\lambda > 0$

- $|\bar{y}_{jk}(\lambda)| := \max\{|y_{j'k'}|; \ I_{j'k'} \in \mathcal{T}_{jk}^{(\eta)}(\lambda)\}$,
- $|\bar{\beta}_{jk}(\lambda)| := \max\{|\beta_{j'k'}|; \ I_{j'k'} \in \mathcal{T}_{jk}^{(\eta)}(\lambda)\}$,
- $|\bar{\delta}_{jk}(\lambda)| := \max\{|y_{j'k'} - \beta_{j'k'}|; \ I_{j'k'} \in \mathcal{T}_{jk}^{(\eta)}(\lambda)\}$.

**Remark 8.1.** *For any $\lambda > 0$,*

$$|\bar{\beta}_{jk}(\lambda)| \leq \frac{\lambda}{2} \iff \forall \ I_{j'k'} \in \mathcal{T}_{jk}^{(\eta)}(\lambda), \ |\beta_{j'k'}| \leq \frac{\lambda}{2},$$

$$|\bar{\beta}_{jk}(\lambda)| > \frac{\lambda}{2} \iff \exists \ I_{j'k'} \in \mathcal{T}_{jk}^{(\eta)}(\lambda), \ |\beta_{j'k'}| > \frac{\lambda}{2}.$$

Note that $|\bar{y}_{jk}(.)|, |\bar{\beta}_{jk}(.)|$ and $|\bar{\delta}_{jk}(.)|$ are decreasing functions with respect to $\lambda$.

Choosing $m^2 \geq 16\eta$, we have

$$\sum_{j=0}^{j_{\lambda_\epsilon}-1} \sum_k \beta_{jk}^2 \mathbf{1}\left\{\forall \ I_{j'k'} \in \mathcal{T}_{jk}^{(\eta)}(\lambda_\epsilon), |\beta_{j'k'}| \leq \frac{\lambda_\epsilon}{2}\right\}$$

$$= \sum_{j=0}^{j_{\lambda_\epsilon}-1} \sum_k \beta_{jk}^2 \mathbf{1}\left\{|\bar{\beta}_{jk}(\lambda_\epsilon)| \leq \frac{\lambda_\epsilon}{2}\right\}$$

$$= \mathbb{E} \sum_{j=0}^{j_{\lambda_\epsilon}-1} \sum_k \beta_{jk}^2 \mathbf{1}\left\{|\bar{\beta}_{jk}(\lambda_\epsilon)| \leq \frac{\lambda_\epsilon}{2}\right\} [\mathbf{1}\{|\bar{y}_{jk}(\lambda_\epsilon)| \leq \lambda_\epsilon\} + \mathbf{1}\{|\bar{y}_{jk}(\lambda_\epsilon)| > \lambda_\epsilon\}]$$

$$
\begin{aligned}
\leq \quad & \mathbb{E} \sum_{j=0}^{j_{\lambda_\epsilon}-1} \sum_k (\beta_{jk} - \gamma_{jk} y_{jk})^2 \mathbf{1}\{|\bar{y}_{jk}(\lambda_\epsilon)| \leq \lambda_\epsilon\} \\
& + \mathbb{E} \sum_{j=0}^{j_{\lambda_\epsilon}-1} \sum_k \beta_{jk}^2 \mathbb{P}(|\bar{y}_{jk}(\lambda_\epsilon)| > \lambda_\epsilon) \mathbf{1}\left\{|\bar{\beta}_{jk}(\lambda_\epsilon)| \leq \frac{\lambda_\epsilon}{2}\right\} \\
\leq \quad & \mathbb{E} \sum_{j=0}^{j_{\lambda_\epsilon}-1} \sum_k (\beta_{jk} - \gamma_{jk} y_{jk})^2 + C\lambda_\epsilon^2 \sum_{j=0}^{j_{\lambda_\epsilon}-1} \sum_k \mathbb{P}(|\bar{\delta}_{jk}(\lambda_\epsilon)| > \frac{\lambda_\epsilon}{2}) \\
\leq \quad & \mathbb{E}\|\tilde{f}_T - f\|_2^2 + C \, j_{\lambda_\epsilon} \, 2^{j_{\lambda_\epsilon}} \lambda_\epsilon^2 \epsilon^{\frac{m^2}{8}} \\
\leq \quad & \mathbb{E}\|\tilde{f}_T - f\|_2^2 + C\lambda_\epsilon^2 \\
\leq \quad & C\lambda_\epsilon^{\frac{4s}{1+2s}}.
\end{aligned}
$$

So, using the continuity of $\lambda_\epsilon$ in 0, we deduce that

$$
\sup_{\lambda > 0} \lambda^{-\frac{4s}{1+2s}} \sum_{j=0}^{j_\lambda - 1} \sum_k \beta_{jk}^2 \mathbf{1}\left\{\forall \, I_{j'k'} \in \mathcal{T}_{jk}^{(\eta)}(\lambda), |\beta_{j'k'}| \leq \frac{\lambda}{2}\right\} \; < \infty.
$$

It comes that $f \in W_{\frac{2}{1+2s}, \eta}^{T}$. So (8.1) is proved.

*STEP 2:* Let $\epsilon_m > 0$ be such that $\epsilon_m \sqrt{\log(\frac{1}{\epsilon_m})} < m^{-1}$. It suffices to prove that for any $0 < \epsilon \leq \epsilon_m$

$$
\mathbb{E}\|\tilde{f}_T - f\|_2^2 = C \, \epsilon^2 + \mathbb{E}\left\|\sum_{j=0}^{j_{\lambda_\epsilon}-1} \sum_k (\gamma_{jk} y_{jk} - \beta_{jk}) \psi_{j,k}\right\|_2^2 + \sum_{j=j_{\lambda_\epsilon}}^{\infty} \sum_k \beta_{jk}^2 \leq C \, \lambda_\epsilon^{\frac{4s}{1+2s}}.
$$

The term $\displaystyle\sum_{j=j_{\lambda_\epsilon}}^{\infty} \sum_k \beta_{jk}^2$ can be bounded by $C\lambda_\epsilon^{\frac{4s}{1+2s}}$, by using the definition of the Besov space $\mathcal{B}_{2,\infty}^{\frac{s}{\eta(1+2s)}}$.

The term $\mathbb{E} \displaystyle\sum_{j=0}^{j_{\lambda_\epsilon}-1} \sum_k (\gamma_{jk} y_{jk} - \beta_{jk})^2$ can be bounded by $C + D$, where

$$
\begin{aligned}
C + D \quad = \quad & \mathbb{E} \sum_{j=0}^{j_{\lambda_\epsilon}-1} \sum_k \beta_{jk}^2 \, \mathbf{1}\{|\bar{y}_{jk}(\lambda_\epsilon)| \leq \lambda_\epsilon\} \\
& + \mathbb{E} \sum_{j=0}^{j_{\lambda_\epsilon}-1} \sum_k (y_{jk} - \beta_{jk})^2 \, \mathbf{1}\{|\bar{y}_{jk}(\lambda_\epsilon)| > \lambda_\epsilon\}.
\end{aligned}
$$

We split $C$ into $C_1 + C_2$ as follows:

$$
\begin{aligned}
C_1 &= \mathbb{E} \sum_{j=0}^{j_{\lambda_\epsilon}-1} \sum_k \beta_{jk}^2 \, \mathbf{1}\{|\bar{y}_{jk}(\lambda_\epsilon)| \leq \lambda_\epsilon\} \mathbf{1}\{|\bar{\beta}_{jk}(\lambda_\epsilon)| \leq 2\lambda_\epsilon\} \\
C_2 &= \mathbb{E} \sum_{j=0}^{j_{\lambda_\epsilon}-1} \sum_k \beta_{jk}^2 \, \mathbf{1}\{|\bar{y}_{jk}(\lambda_\epsilon)| \leq \lambda_\epsilon\} \mathbf{1}\{|\bar{\beta}_{jk}(\lambda_\epsilon)| > 2\lambda_\epsilon\}.
\end{aligned}
$$

Since $f \in W_{\frac{2}{1+2s},\eta}^{T}$ and $f \in \mathcal{B}_{2,\infty}^{\frac{s}{\eta(1+2s)}}$,

$$
\begin{aligned}
C_1 &= \mathbb{E} \sum_{j=0}^{j_{\lambda_\epsilon}-1} \sum_k \beta_{jk}^2 \, \mathbf{1}\{|\bar{y}_{jk}(\lambda_\epsilon)| \leq \lambda_\epsilon\} \mathbf{1}\{|\bar{\beta}_{jk}(\lambda_\epsilon)| \leq 2\lambda_\epsilon\} \\
&\leq \sum_{j=0}^{j_{\lambda_\epsilon}-5-\lfloor \log_2(\eta)\rfloor} \sum_k \beta_{jk}^2 \, \mathbf{1}\{|\bar{\beta}_{jk}(4\lambda_\epsilon)| \leq 2\lambda_\epsilon\} + \sum_{j=j_{\lambda_\epsilon}-4-\lfloor \log_2(\eta)\rfloor}^{\infty} \sum_k \beta_{jk}^2 \\
&\leq C \lambda_\epsilon^{\frac{4s}{1+2s}}
\end{aligned}
$$

and

$$
\begin{aligned}
C_2 &= \mathbb{E} \sum_{j=0}^{j_{\lambda_\epsilon}-1} \sum_k \beta_{jk}^2 \, \mathbf{1}\{|\bar{y}_{jk}(\lambda_\epsilon)| \leq \lambda_\epsilon\} \mathbf{1}\{|\bar{\beta}_{jk}(\lambda_\epsilon)| > 2\lambda_\epsilon\} \\
&\leq \sum_{j=0}^{j_{\lambda_\epsilon}-1} \sum_k \beta_{jk}^2 \, \mathbb{P}(|\bar{\delta}_{jk}(\lambda_\epsilon)| > \lambda_\epsilon) \\
&\leq C \, j_{\lambda_\epsilon} \, 2^{j_{\lambda_\epsilon}} \epsilon^{\frac{m^2}{2}} \\
&\leq C \lambda_\epsilon^{\frac{4s}{1+2s}}.
\end{aligned}
$$

We have used here the concentration property of the Gaussian distribution and the fact that $m^2 \geq 4(1+\eta)$.

We split $D$ into $D_1 + D_2$ as follows:

$$
\begin{aligned}
D_1 &= \mathbb{E} \sum_{j=0}^{j_{\lambda_\epsilon}-1} \sum_k (y_{jk}-\beta_{jk})^2 \, \mathbf{1}\{|\bar{y}_{jk}(\lambda_\epsilon)| > \lambda_\epsilon\} \mathbf{1}\left\{|\bar{\beta}_{jk}(\lambda_\epsilon)| \leq \frac{\lambda_\epsilon}{2}\right\} \\
D_2 &= \mathbb{E} \sum_{j=0}^{j_{\lambda_\epsilon}-1} \sum_k (y_{jk}-\beta_{jk})^2 \, \mathbf{1}\{|\bar{y}_{jk}(\lambda_\epsilon)| > \lambda_\epsilon\} \mathbf{1}\left\{|\bar{\beta}_{jk}(\lambda_\epsilon)| > \frac{\lambda_\epsilon}{2}\right\}.
\end{aligned}
$$

For $D_1$ we use the Cauchy-Schwartz inequality:

$$
\mathbb{E}(y_{jk}-\beta_{jk})^2 \mathbf{1}\left\{|\bar{\delta}_{jk}(\lambda_\epsilon)| > \frac{\lambda_\epsilon}{2}\right\} \leq 2^{\frac{j_{\lambda_\epsilon}}{2}} (\mathbb{P}\left(|y_{jk}-\beta_{jk}| > \frac{\lambda_\epsilon}{2}\right)^{1/2} (\mathbb{E}(y_{jk}-\beta_{jk})^4)^{1/2}
$$

where $\mathbb{E}(y_{jk} - \beta_{jk})^4 = 3\epsilon^4$ and $\mathbb{P}(|y_{jk} - \beta_{jk}| > \frac{\lambda_\epsilon}{2}) \leq \epsilon^{\frac{m^2}{8}}$ (using the concentration properties of the Gaussian distribution). So, choosing $m$ such that $m^2 \geq 48\eta$,

$$
\begin{aligned}
D_1 &\leq C \; 2^{\frac{j_{\lambda_\epsilon}}{2}} \sum_{j=0}^{j_{\lambda_\epsilon}-1} \sum_k \epsilon^2 \mathbf{1}\left\{ |\bar{\beta}_{jk}(\lambda_\epsilon)| \leq \frac{\lambda_\epsilon}{2} \right\} \epsilon^{\frac{m^2}{16}} \\
&\leq C \; 2^{\frac{3j_{\lambda_\epsilon}}{2}} \lambda_\epsilon^{2 + \frac{m^2}{16}} \\
&\leq C \lambda_\epsilon^{\frac{4s}{1+2s}}.
\end{aligned}
$$

For $D_2$, we use Proposition 5.2 with $r = \frac{2}{1+2s}$ and $\lambda_0 = m\epsilon_m \sqrt{\log(\frac{1}{\epsilon_m})}$:

$$
\begin{aligned}
D_2 &\leq \sum_{j=0}^{j_{\lambda_\epsilon}-1} \sum_k \epsilon^2 \mathbf{1}\left\{ |\bar{\beta}_{jk}(\lambda_\epsilon)| > \frac{\lambda_\epsilon}{2} \right\} \\
&\leq \sum_{j=0}^{j_{\lambda_\epsilon}-1} \sum_k \epsilon^2 \mathbf{1}\left\{ \exists I_{j'k'} \in \mathcal{T}_{jk}^{(\eta)}(\lambda_\epsilon) \; / \; |\beta_{j'k'}| > \frac{\lambda_\epsilon}{2} \right\} \\
&\leq C \lambda_\epsilon^{\frac{4s}{1+2s}}.
\end{aligned}
$$

Looking at bounds of $C_1$, $C_2$, $D_1$ and $D_2$, (8.2) is proved. $\qquad\square$